%% file: McRobieGPD14April2013arxiv1.tex
\documentclass[a4paper]{article}
\usepackage{natbib}

  \usepackage{graphics,graphicx}
   \usepackage{amsmath}
   \usepackage{amsfonts}
   \usepackage{float,txfonts,amssymb,bm}
 \usepackage{colortbl}

\newcommand{\Frechet}{Frech\'{e}t }

\newcommand{\hs}{\hspace*{1.5cm}}
\newcommand{\hsf}{\hspace*{0.47cm}}
\newcommand{\hsff}{\hspace*{0.8cm}}
\newcommand{\hsfg}{\hspace*{0.65cm}}

\newcommand{\cs}{\cellcolor[gray]{0.8}}

 \bmdefine\bone{1}    \bmdefine\bX{X}  \bmdefine\bt{t}
\bmdefine\btau{\tau} \bmdefine\bx{x}  \bmdefine\balpha{\alpha}
\bmdefine\bbeta{\beta} \bmdefine\by{y} \bmdefine\bu{u}
\bmdefine\bG{G}

\begin{document}

\title{Elemental unbiased estimators for the Generalized Pareto tail}
\author{Allan McRobie \\
Cambridge University Engineering Department\\
Trumpington St, Cambridge, CB2 1PZ, UK \\
fam20@cam.ac.uk}

\maketitle

\begin{abstract}
Unbiased location- and scale-invariant `elemental' estimators for the GPD tail
parameter are constructed.
Each involves three log-spacings.
The estimators are unbiased for finite sample sizes, even as small as $N=3$.
It is shown that the elementals form a complete basis for unbiased location- and scale-invariant estimators
constructed from linear combinations of log-spacings.
Preliminary numerical evidence is presented which suggests that elemental combinations can be constructed
which are consistent estimators of the tail parameter for samples drawn from the pure GPD family.
\end{abstract}

\maketitle

\section{Introduction}
The Generalized Pareto Distribution (GPD) and the Generalized
Extreme Value (GEV) distribution play a central role in extreme
value theory. Each has three parameters $(\mu, \sigma, \xi)$
corresponding to location, scale and tail (or shape) respectively.
This paper describes a particularly simple set of location- and scale-invariant
`elemental' estimators for the GPD tail
parameter. Each `elemental' involves three log-spacings of the data, and each
is unbiased over all tail parameters $-\infty < \xi < \infty$, and for all sample sizes, as small as $N=3$.

\begin{figure}[h!]
\begin{picture}(150,40)
\thinlines \multiput(30,30)(10,0){11}{\circle*{1.5}} \put (20,30)
{\line(1,0){120}} \put (30,33){$X_n$} \put(50,33){$X_J$}
\put(60,33){$X_{J-1}$} \put(90,33){$X_{I+1}$} \put(100,33){$X_I$}
\put (120,33){$X_2$} \put(130,33){$X_1$ (Largest)}
\put(50,5){\line(0,1){25}}\put(60,5){\line(0,1){25}}
\put(90,5){\line(0,1){25}}\put(100,5){\line(0,1){25}} \thicklines
\put(50,15){\line(1,0){50}}\thicklines\put(60,20){\line(1,0){40}}
\thicklines\put(50,10){\line(1,0){40}}
\put(50,15){\circle*{2}}\put(100,15){\circle*{2}}
\put(60,20){\circle*{2}}\put(100,20){\circle*{2}}
\put(50,10){\circle*{2}}\put(90,10){\circle*{2}} \put(105,19){$J-1$}
\put(105,14){$-(J-1)+I$} \put(105,9){$-I$} \put(75,11){$t$}
\put(75,17){$\tau$}
\end{picture}
\caption{Between any two non-adjacent data points $X_I$ and $X_J$ an
elemental estimator $\hat{\xi}_{IJ}$ can be defined. It involves
three log-spacings - the one between $X_I$ and $X_J$, together with
two shorter log-spacings connecting each end-point to the data point
immediately inside the other end.} \label{bardiag}
\end{figure}
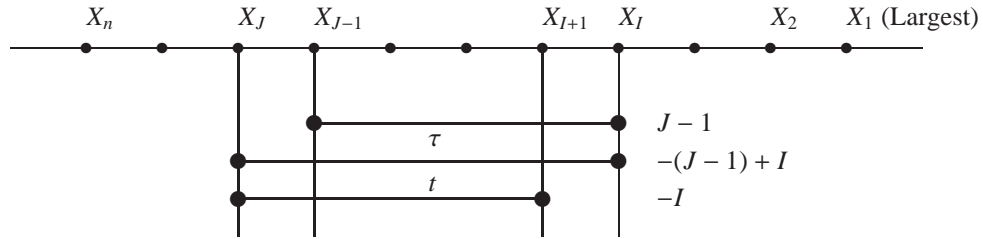

The elemental estimators (illustrated in Figure \ref{bardiag}) have
the form
\begin{equation}
\hat{\xi}_{IJ} =
\log \frac{\tau^{J-1}}{t^I} \ \ \ \text{where} \ \ \ \tau =
\frac{X_I-X_{J-1}}{X_I-X_J} \ \ \ \text{and} \ \ \ t = \frac{X_{I+1}
-X_J}{X_I- X_J}, \ \ \text{with} \  J \geq I+2 \label{equation1}
\end{equation}
and the $X_I$ are the upper-order statistics, numbered in decreasing
order starting from $I=1$ as the data maximum.

The Generalized Pareto Distribution arises as the limiting
distribution of maxima in Peaks-Over-Threshold approaches (see for
example \citet{embrechts}). It has distribution function:
\begin{equation}
F(x) = 1-\left( 1+ \xi \frac{x-\mu}{\sigma} \right)^{-1/\xi}
\label{thegpd}
\end{equation}
The parameters $\mu$ and $\xi$ can take any value on the real line,
whilst $\sigma$ can be any positive value (and when $\xi =0$ the
distribution function (\ref{thegpd}) becomes the exponential
distribution). For GPDs with positive $\xi$, the support ($\mu \leq x$) is
unbounded at the right, giving long- or heavy-tailed distributions.
For $\xi$ negative, the support is bounded both below and above ($ \mu \leq x \leq \mu - \sigma / \xi$).

\section{Other estimators}
Estimators for the tail parameter can be loosely classed into:
maximum likelihood (ML); method of moments; Pickands-like and
Bayesian. Standard texts such as \cite{embrechts} and \cite{reiss}
provide detailed background, with \cite{coles} giving the Bayesian
perspective. \cite{bermudez1} provide a comprehensive review, such
that only a brief survey is presented here.

The maximum likelihood approach to the GPD is described in
\cite{smith}. Although it possesses some desirable properties, the
numerical maximization algorithms can experience problems for small
sample sizes and for negative tail parameters, as there are
occasions when the likelihood function does not possess a local
maximum (\cite{castillo}). To avoid such problems, a method of
moments approach was proposed by \cite{hosking}.

The classical tail parameter estimator is that of \cite{hill}.
However, it is not location invariant and is only valid in the
heavy-tailed \Frechet region ($\xi$ positive) of the GEV, although
an extension into the Weibull region ($\xi$ negative) was proposed
by \cite{Dekkers} using the method of moments. \cite{pickands}
proposed an estimator based on log-spacings which overcame many of
these shortcomings. This estimator is popular in current
applications, and a substantial literature exists on its
generalization (\cite{drees2, yun}, for example), the most general
and efficient of which appear to be those of \cite{segers}, derived
using second order theory of regular variation. These are optimised
for estimation of the tail index in the more general case of data
drawn from {\it any} distribution within the domain of attraction of
the particular GPD. 
Although the main concern of Extreme Value Theory is the domain of attraction case,
this paper restricts attention to distributions within the pure GPD family.
The possibility that results derived in this specific setting
may be extended to the more general case is left for later consideration.

Throughout, there is an emphasis on results that are valid for small sample sizes.

\section{Elemental Estimators}

\begin{figure}[h!] \centering

  \includegraphics[width=75mm,keepaspectratio]{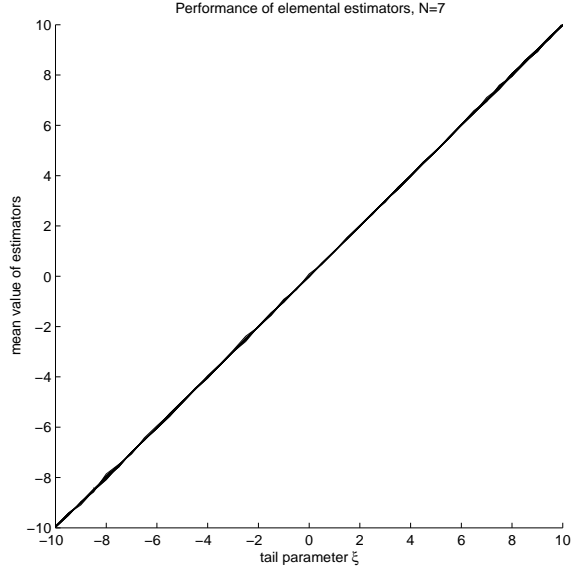}
  \caption{Averages over 50,000 samples for the 15 elemental estimators
  for $N=7$ over a range of $\xi$. The fifteen lines are almost indistinguishable
  from the diagonal, indicating that each is indeed absolutely unbiased.
  }\label{pickands1}
\end{figure}

The main result here is the proof in Appendix 1 that each elemental
estimator is absolutely unbiased within the GPD family. 
That the proof is valid, remarkably, for ALL $\xi$ may be
appreciated by inspection of Eqn.~\ref{squigneg} there. For $\xi$
negative, the expectation of the log-spacing is expressed in terms
of the tail probabilities $G_i$ and $G_j$ via a term
$\log(G_j^\gamma - G_i^\gamma)$ with $\gamma = -\xi$. This trivially
decomposes into a simple term $\gamma \log G_j$ and a complicated
term $\log(1-(G_i/G_j)^\gamma)$. The proof shows how the simple
terms provide the expectation $\hat{\xi} = - \gamma$, and how the
elementals combine the complicated terms in such a manner that they
cancel (obviating the need to evaluate them explicitly). For $\xi$
positive, the absolute lack of bias is maintained by an additional
simple term $\gamma \log G_i G_j$ which adds $2\gamma$ to the
$-\gamma$ result for $\xi$ negative. This elegant correspondence
between the results for $\xi$ positive and negative is absent in
previous approaches to GPD tail estimation. 
As further demonstration of the absolute lack of bias of each elemental triplet, even at small
sample sizes, the numerically-returned average values for each of
the fifteen elemental estimators available for $N= 7$ are shown in
Fig.~\ref{pickands1} over a wide range of tail parameters ($-10 \leq
\xi \leq 10$).

\vspace*{1mm}

\section{Linear combinations}

Trivially, any unit-weight linear combination of elementals will
also be unbiased. Whilst it will be of interest to form efficient
combinations, no detailed analysis of efficiency or variance is
undertaken here. Instead, the performance of a simple combination is
reported.

Linear combinations of elementals are most conveniently described
via the upper triangular matrix $M$ of dimension $N\times N$ of
Table \ref{delmatrix} containing all possible log-spacings. The
general term is $M_{IJ} = \log(X_I-X_{J})$ for $J \geq I+1$ and zero
otherwise. (Use of the $N\times N$ form with zero diagonal allows
for easier indexing).  Each element above the secondary diagonal ($J
\geq I+2$) may be uniquely identified with an elemental estimator,
involving that log-spacing, that to its left and that below it in
$M$. Corresponding weights of elemental estimator combinations may
thus be stored in an $N\times N$ upper triangular matrix $R$.
Weighting an elemental $\hat{\xi}_{IJ}$ by $r_{IJ}$ requires (from
Eqn. \ref{equation1}) that the three weights $r_{IJ}\times \lbrace J-1,
-(J-1-I), -I \rbrace $ be given, respectively,  to the corresponding left, upper-right and
lower log-spacings in $M$. These weights are illustrated in the grid
$G$ shown in Table \ref{picktable2}. The totals may then be collected in
an $N \times N$ matrix $A$ of log-spacing weights. The sum of the
entrywise product of $A$ and $M$ is then the unbiased estimate
$\hat{\xi}$.

\begin{table}[!]

\vspace*{5mm}

\begin{tabular}{|c|c|c|c|c|c|c|}
  \hline
\hspace{1.5cm} & $\log (X_1-X_2) $ & $\log (X_1-X_3)$   & $\log
(X_1-X_4)$ &
$\log (X_1-X_5)$   & $\log(X_1-X_6)$   & $\log (X_1-X_7)$\\
\hline
 &         & $\log (X_2-X_3)$ & $\log(X_2-X_4)$ & $\log(X_2-X_5)$ & $\log (X_2-X_6)$ & $\log(X_2-X_7)$ \\ \hline
  &        &           & $\log(X_3-X_4)$ & $ \cs \log(X_3-X_5)$ & \cs $\log(X_3-X_6)$ & $\log(X_3-X_7)$ \\ \hline
  &        &           &           & $\log(X_4-X_5)$ & \cs $\log(X_4-X_6)$ & $\log(X_4-X_7)$ \\ \hline
  &        &           &           &           & $\log(X_5-X_6)$ & $\log(X_5-X_7)$ \\ \hline
  &        &           &           &           &           & $\log(X_6-X_7)$ \\ \hline
  &        &           &           &           &           & \\ \hline
\end{tabular}
\caption{The address matrix $M$ for $N=7$. Each elemental involves
three adjacent cells in an inverted-L formation. That corresponding
to
 the elemental $\hat{\xi}_{36}$ is shaded for illustration. \label{delmatrix}}

\vspace*{10mm}

  \begin{tabular}{|c|lcr|lcr|lcr|lcr|lcr|lcr|}
  \hline
 \hs & \hsff &  & 2 & -1 & \hsf & 3 & -2 & \hsf & 4 & -3 & \hsf & 5 & -4 & \hsf & 6 & -5 & \hsfg &  \\ \hline
 \hs & \hsf &  &   & -1 & &   & -1 & &   & -1 & &   & -1 & &   & -1 &  &  \\
 \hs & \hsf &  &   &    & &   &    & &   &    & &   &    & &   &    &  &  \\
 \hs & \hsf &  &   &    & & 3 & -1 & & 4 & -2 & & 5 & -3 & & 6 & -4 &  &  \\ \hline
 \hs & \hsf &  &   &    & &   & -2 & &   & -2 & &   & -2 & &   & -2 &  &  \\
 \hs & \hsf &  &   &    & &   &    & &   &    & &   &    & &   &    &  &  \\
 \hs & \hsf &  &   &    & &   &    & & 4 & -1 & & \cs 5 & \cs -2 & & 6 & -3 &  &  \\ \hline
 \hs & \hsf &  &   &    & &   &    & &   & -3 & &   & \cs
  -3 & &   & -3 &  &  \\
 \hs & \hsf &  &   &    & &   &    & &   &    & &   &    & &   &    &  &  \\
 \hs & \hsf &  &   &    & &   &    & &   &    & & 5 & -1 & & 6 & -2 &  &  \\ \hline
 \hs & \hsf &  &   &    & &   &    & &   &    & &   & -4 & &   & -4 &  &  \\
 \hs & \hsf &  &   &    & &   &    & &   &    & &   &    & &   &    &  &  \\
 \hs & \hsf &  &   &    & &   &    & &   &    & &   &    & & 6 & -1 &  &  \\ \hline
 \hs & \hsf &  &   &    & &   &    & &   &    & &   &    & &   & -5 &  &  \\  \hline
 \hs & \hsf &  &   &    & &   &    & &   &    & &   &    & &   &    &  &  \\ \hline
\end{tabular}

\caption{ The grid $G$ of exponents of the elemental estimators, for
$N=7$. Again, the terms involved in
 the elemental $\hat{\xi}_{36}$ are shaded for illustration.\label{picktable2}}
 \vspace*{5mm}
\end{table}

\vspace*{1mm}

In summary, the matrix $M$ is the roadmap of all log-spacings and
the grid $G$ gives the set of weights to be used within each
elemental. A linear combination of elementals is then defined by a
unit-sum matrix $R$, and the corresponding log-spacing weights are
collected in the zero-sum matrix $A$.

\subsection*{An example combination}

A natural choice of linear combination might give equal weight to
each elemental. However here we give further consideration to a
simple ``linearly-rising'' combination. 
 Numerical
experiments indicate that many simple choices of linear combination
lead to good estimators, and further research may seek the optimal
combination. There is thus nothing special about the
``linearly-rising'' combination considered here. As will be
demonstrated, it has a good all-round performance, but more
importantly it illustrates the great simplicity that the elementals
permit, allowing the ready creation of unbiased tail estimators with
efficiencies comparable to current leading (and often highly
complicated) alternatives.

\begin{figure}[h!] \centering
  \includegraphics[width=45mm,keepaspectratio]{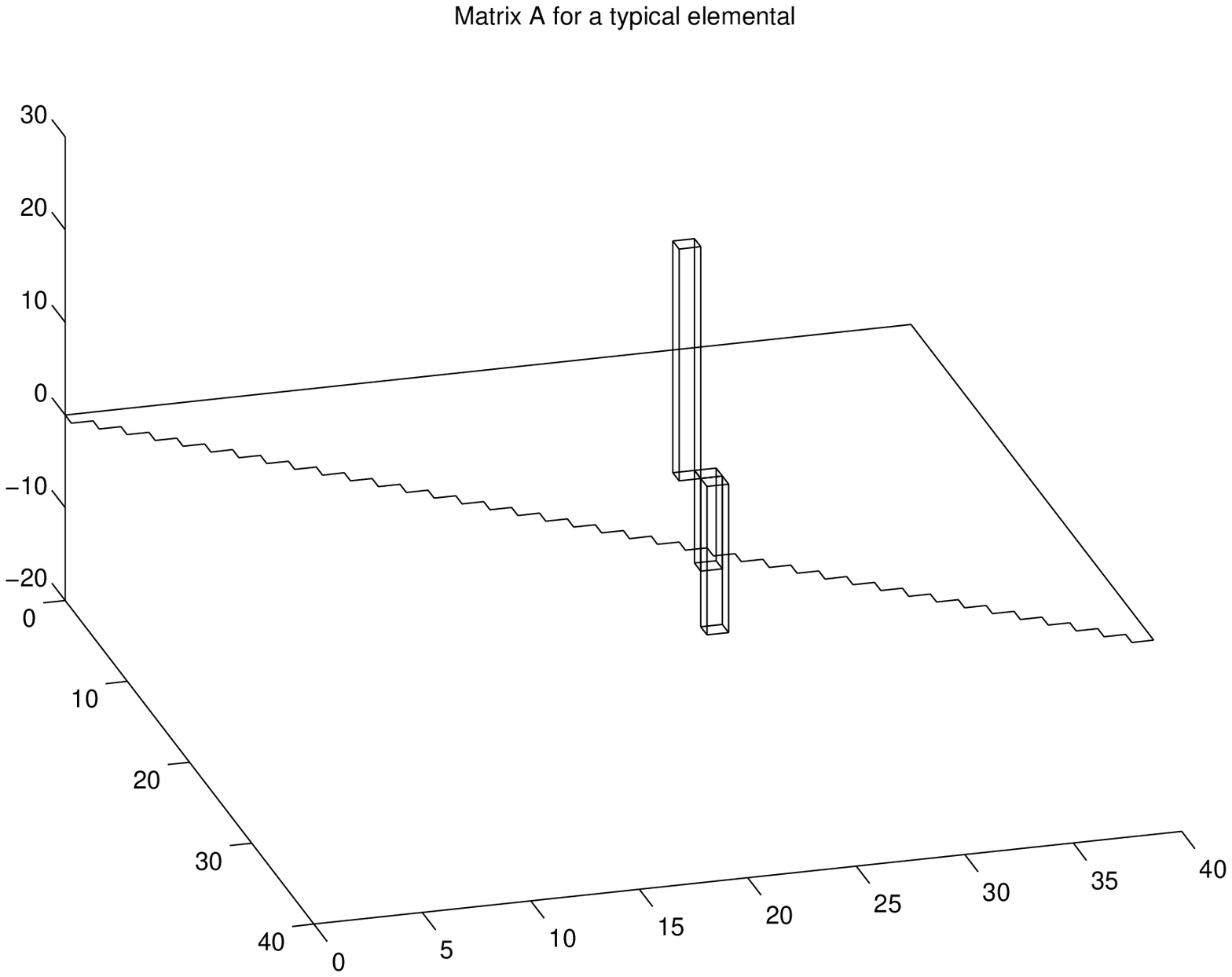}
  \includegraphics[width=45mm,keepaspectratio]{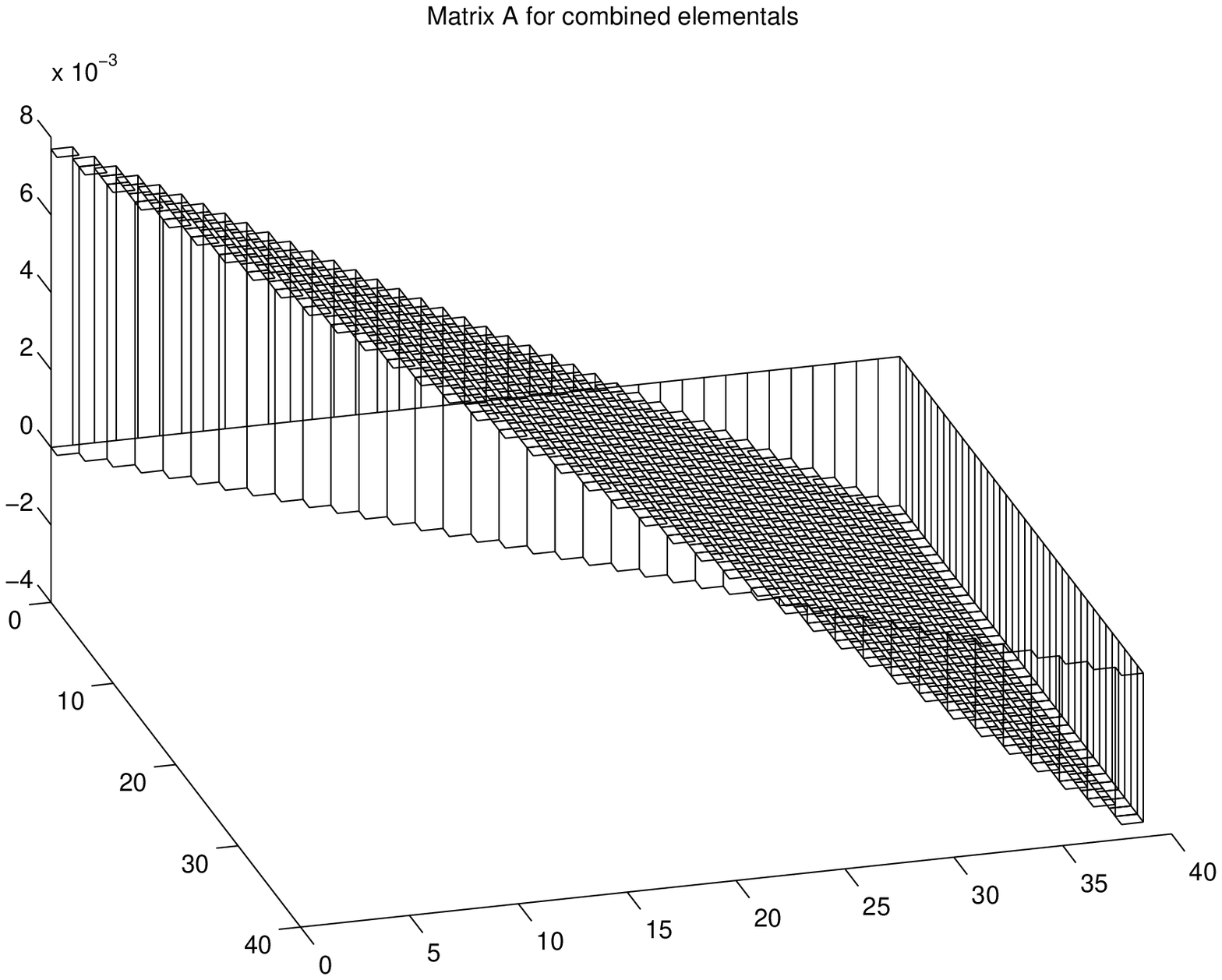}
  \includegraphics[width=45mm,keepaspectratio]{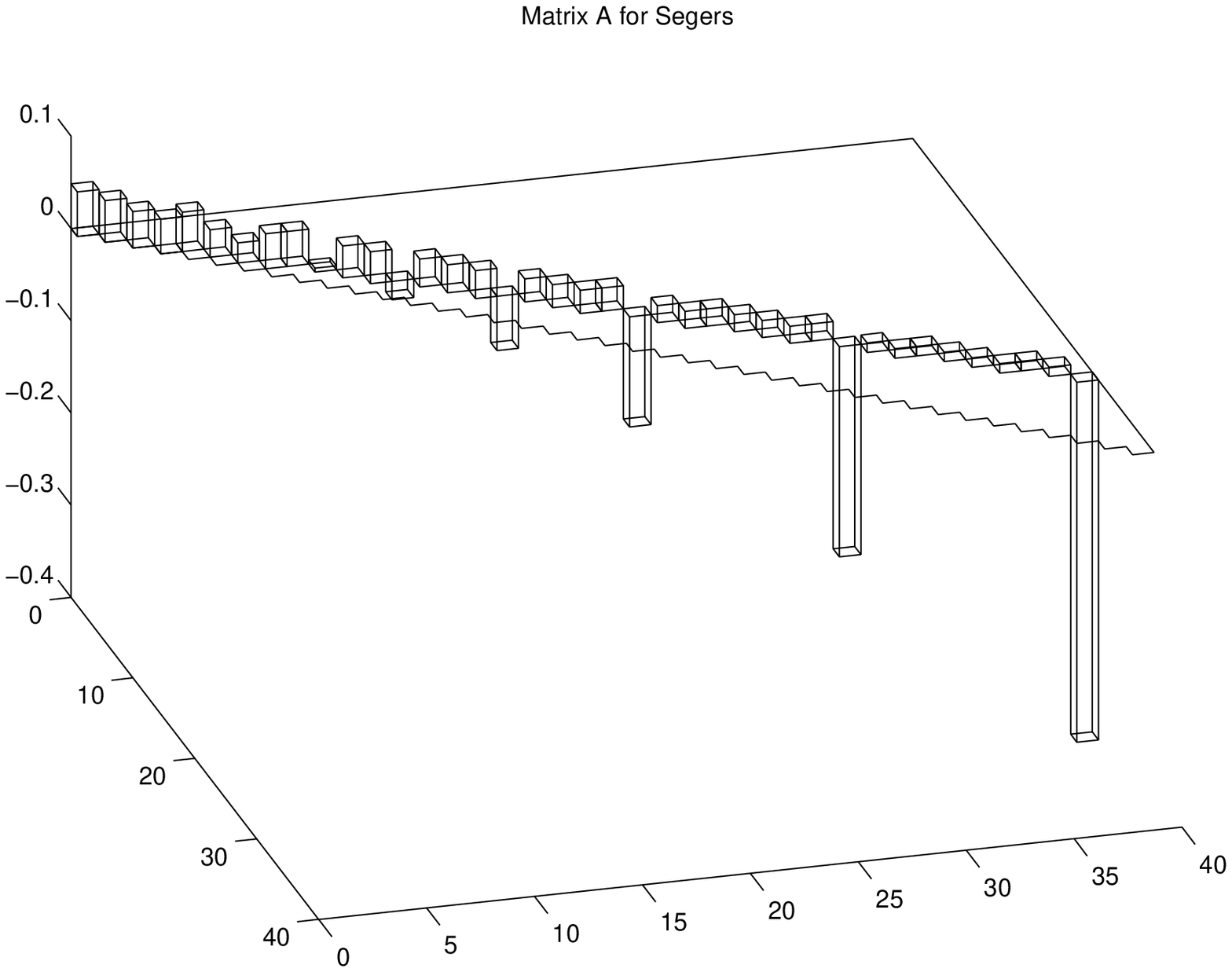}
  \caption{The $A$ matrix of log-spacing weights for a typical elemental,
  the combined elementals and a Segers estimator.
  }\label{Amatrix}
\end{figure}

The ``linearly-rising'' combination has elemental weights $r_{IJ}
\propto N+1-J$ and the resulting log-spacing weights are $a_{IJ} =
6(2N-3J+2)/(N(N-1)(N-2))$ for $J \geq I+1$. For example, for $N=7$
the weights are
\begin{equation}
  R = \  \frac{1}{35} \left(
  \begin{array}{ccccccc}
  . & . & 5 & 4 & 3 & 2 & 1 \\
   & &   & 4 & 3 & 2 & 1 \\
   & &   &   & 3 & 2 & 1 \\
   & &   &   &   & 2 & 1 \\
   & &   &   &   &   & 1 \\
   & &   &   &   &   &   \\
   & &   &   &   &   &   \\
  \end{array}
\right) \ \ \text{giving} \
A = \ \frac{1}{35} \left(
  \begin{array}{ccccccc}
. &   10 & 7 & 4 & 1 & -2 & -5 \\
 &      & 7 & 4 & 1 & -2 & -5 \\
 &      &   & 4 & 1 & -2 & -5 \\
 &      &   &   & 1 & -2 & -5 \\
 &      &   &   &   & -2 & -5 \\
 &      &   &   &   &    & -5 \\
 &      &   &   &   &    &    \\
  \end{array}
\right) \nonumber \label{pickandsform}
\end{equation}
Further illustration is given in Figure~\ref{Amatrix} (centre) for a
sample size of 40, showing how the log-spacing weights have a simple
linearly-rising distribution with zero mean. For comparison, the
unusual pattern of the corresponding log-spacing weights of an
unconstrained (i.e. optimised but biased) Segers estimator are also shown. Since the Segers
estimator has only a single non-zero weight in each column it
immediately follows that it cannot be constructed from the elemental
triplets.

\begin{figure}[h!] \centering
  \includegraphics[width=65mm,keepaspectratio]{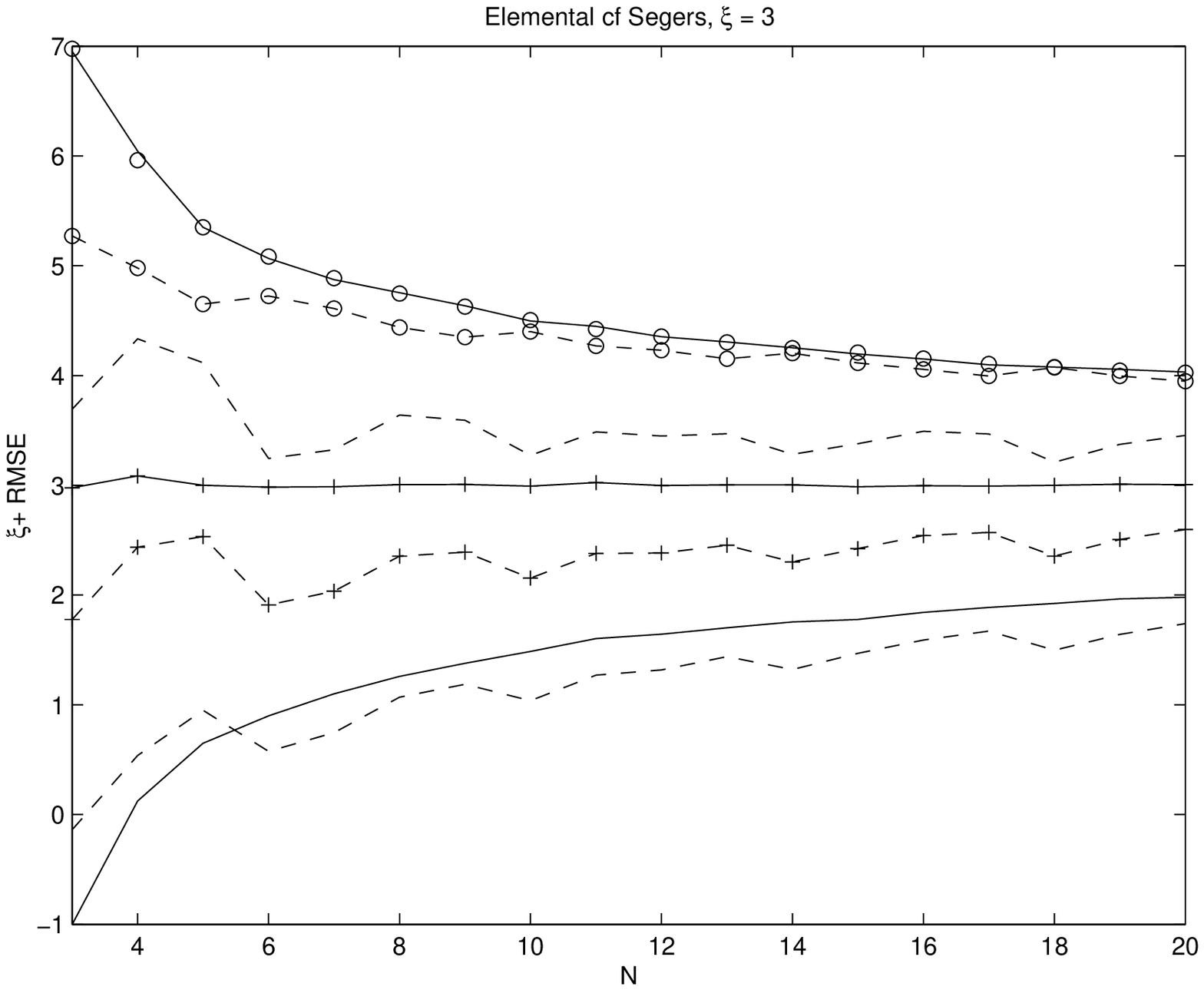}
  \includegraphics[width=65mm,keepaspectratio]{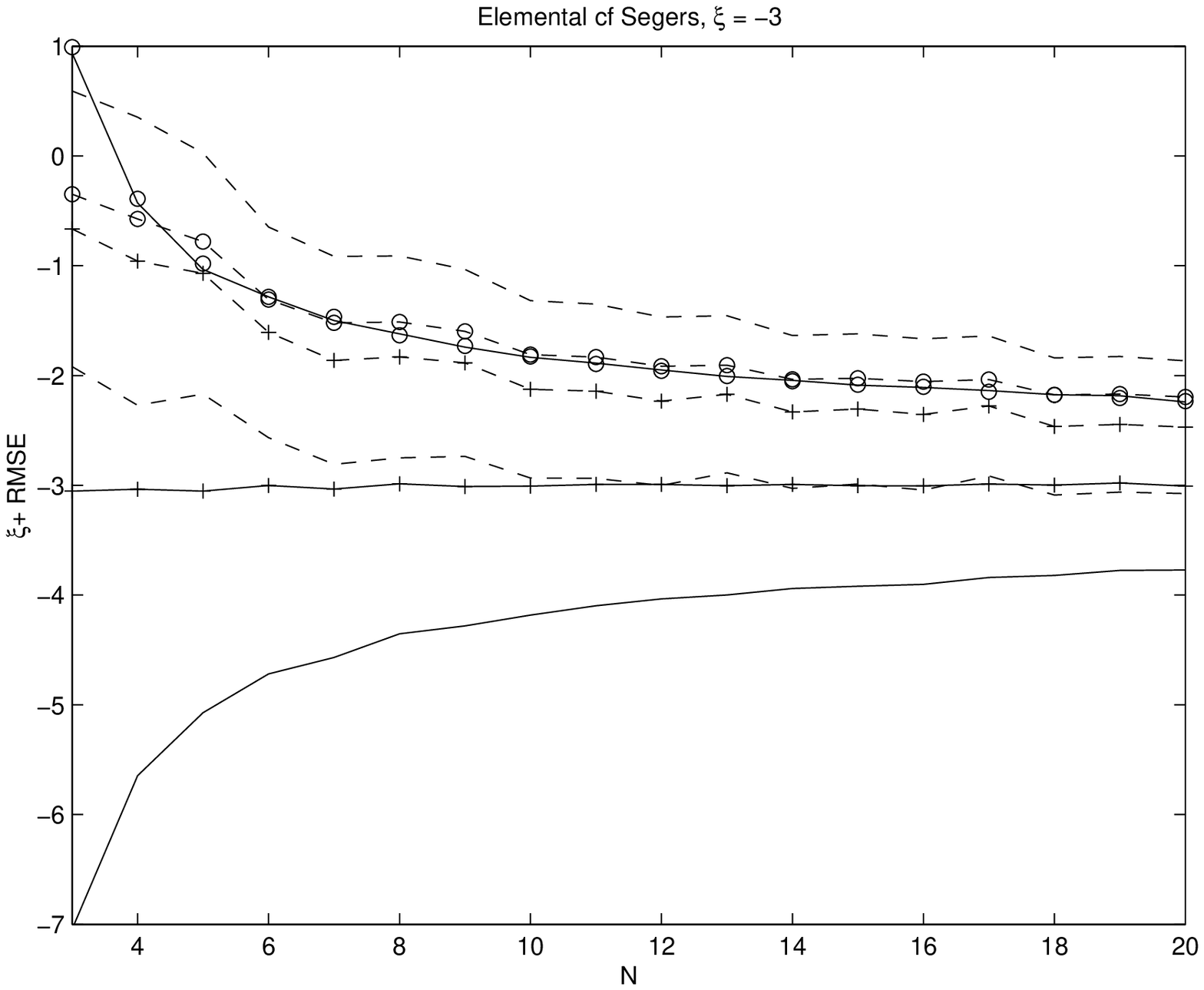}
  \caption{The mean, (mean $\pm$ std. dev.), and (actual + rmse)
  for the elemental combination (solid lines) and the unconstrained Segers (dashed)
  for 10,000 samples of sizes $N = 3$ to $20$ drawn from GPDs with $\xi = 3$ and $\xi = -3$.
  Means are highlighted with $+$ and actual-plus-rmse by $\circ$.
  Note the large bias in the Segers estimator (particularly for $\xi = -3$), and
  note the lack of bias in the elemental combination
  for samples as small as $N=3$.
  }\label{squigpm}
\end{figure}

In Figure~\ref{squigpm}, the errors of the elemental combination are
compared with those of the unconstrained Segers estimator for pure
GPDs with the somewhat extreme cases of $\xi = \pm 3$, and small
sample sizes. The elemental combination has comparatively large
variance around an unbiased mean, in contrast to the Segers
estimator which is more tightly bunched around a biased offset.
Despite the complexity, the extensive optimization and the
substantial bias in the Segers estimator, its mean square error is
nevertheless typically only marginally less than that of the simple
elemental combination.

\section{Completeness}
{\bf Proposition}: if $\hat{\xi}$ is a linear combination of
log-spacings and is an absolutely-unbiased, location- and
scale-invariant estimator of the tail parameter of the GPD, then
$\hat{\xi}$ is a linear combination of elementals.

The proof, presented in Appendix 2, shows that a requirement for
lack of bias imposes $N-1$ independent constraints on the $N(N-1)/2$
dimensional space of possible linear combinations. The resulting
subspace of unbiased estimators thus has dimension $(N-1)(N-2)/2$
and is that subspace spanned by the elementals. The elementals thus
form a complete basis for unbiased, location- and scale-invariant
log-spacing estimators of the GPD tail parameter.

\section{Efficiency and Optimality}
Given that the efficiency of the estimator depends on the actual
unknown value of the parameter $\xi$, there is no unique definition
of optimality. The question as to which linear combination is in
some sense the `best' is thus a matter of judgement.

\begin{figure}[h!] \centering
  \includegraphics[width=95mm,keepaspectratio]{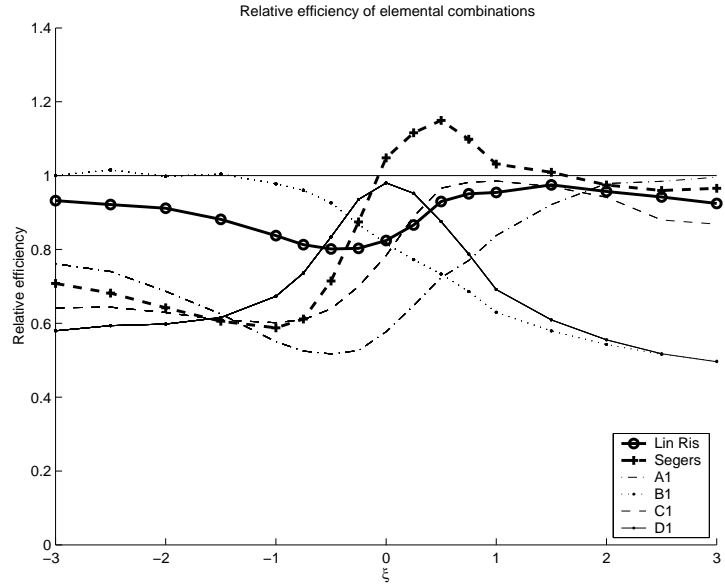}
  \caption{The relative efficiency of various linear combinations of
  elementals (for GPD samples of size $N=20$ over a range of $\xi$).
  Efficiency is defined relative to the numerically-computed minimal
  variance at given $\xi$ within the class of location- and
  scale-invariant unbiased estimators that are linear combinations of log-spacings.
  The linearly-rising combination $r_{ij} \propto N+1-j$ is seen to
  be a good compromise, giving high relative efficiency over the
  whole range $-3 \leq \xi \leq 3$.
  }\label{optimal}
\end{figure}

Fig.~\ref{optimal} shows the relative efficiency of various
linear combinations of elementals,
wherein relative efficiency is defined with
respect to the minimal possible variance (given the tail parameter
$\xi$) within the class of location- and scale-invariant unbiased
estimators which are linear combinations of log-spacings of GPD
data. Using the completeness of the elementals, at any $\xi$ and for
any sample size $N$, the unbiased combination giving minimum
variance within this class can be estimated numerically by
constructing, via repeated samples, the numerical covariance matrix
for the set of $(N-1)(N-2)/2$ elementals, and applying a Lagrange
multiplier to enforce the unit sum condition on the coefficients
$r_{ij}$. The Lagrange multiplier is then an estimate of the minimum
variance. Since the computed coefficients are minimal for that set
of samples, it will, for that sample set, perform better than the
actual global optimum, and thus provide a lower bound on the minimum variance.
The computed coefficients will not be fully optimal for other randomly drawn
sample sets, and since the global optimum is, on average, optimal
for other sets, then an upper bound on the minimum possible variance
(within this class of estimators) can be obtained by applying the
numerically-computed optimal coefficients to a large set of samples
which were not used in their computation. By this procedure, using
two separate blocks of 8000 samples of size $N=20$ drawn from GPDs,
the (approximate) optimal linear combination within the class was constructed for
various $\xi$.

\begin{figure}[h!] \centering
  \includegraphics[width=65mm,keepaspectratio]{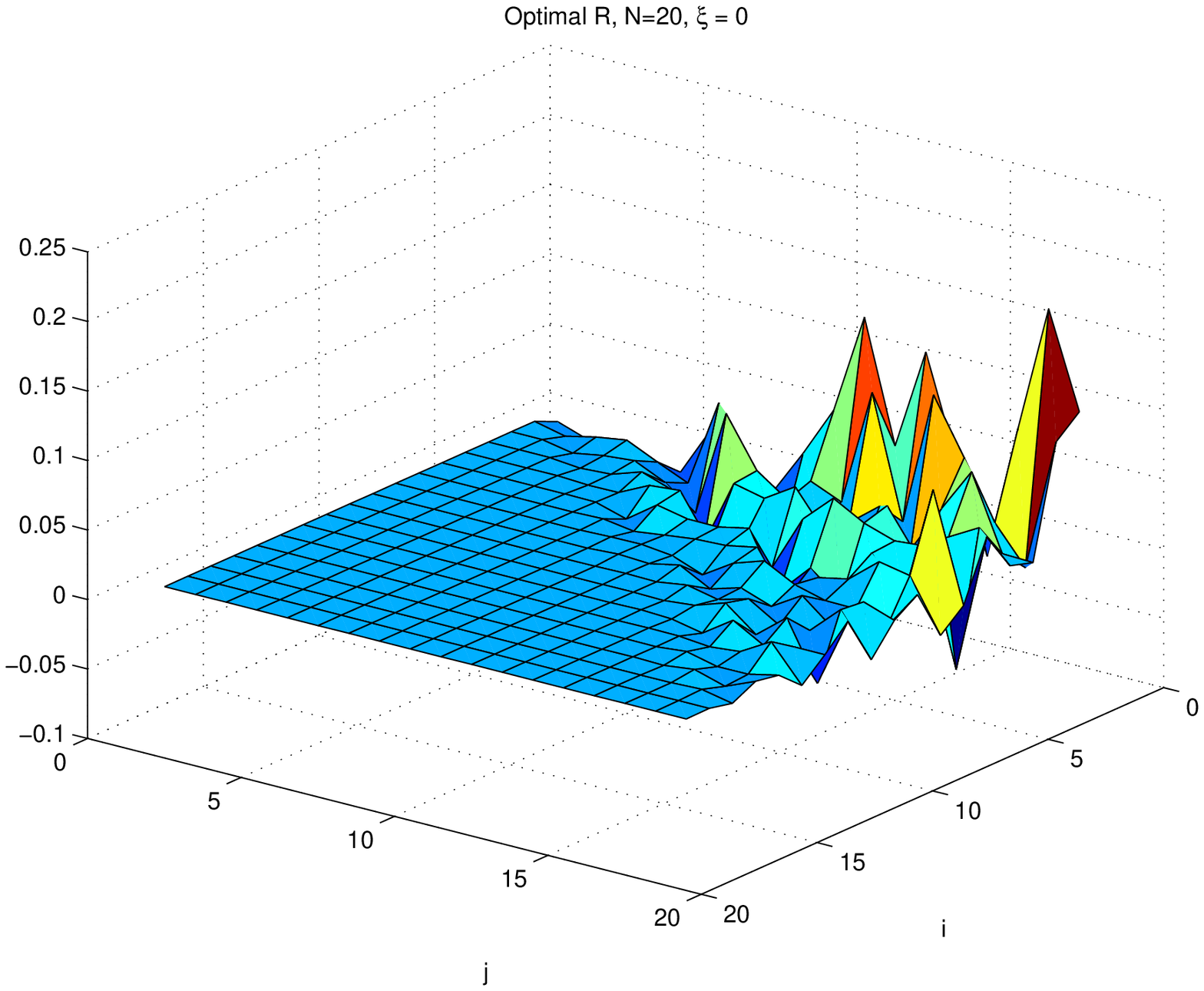}
  \includegraphics[width=65mm,keepaspectratio]{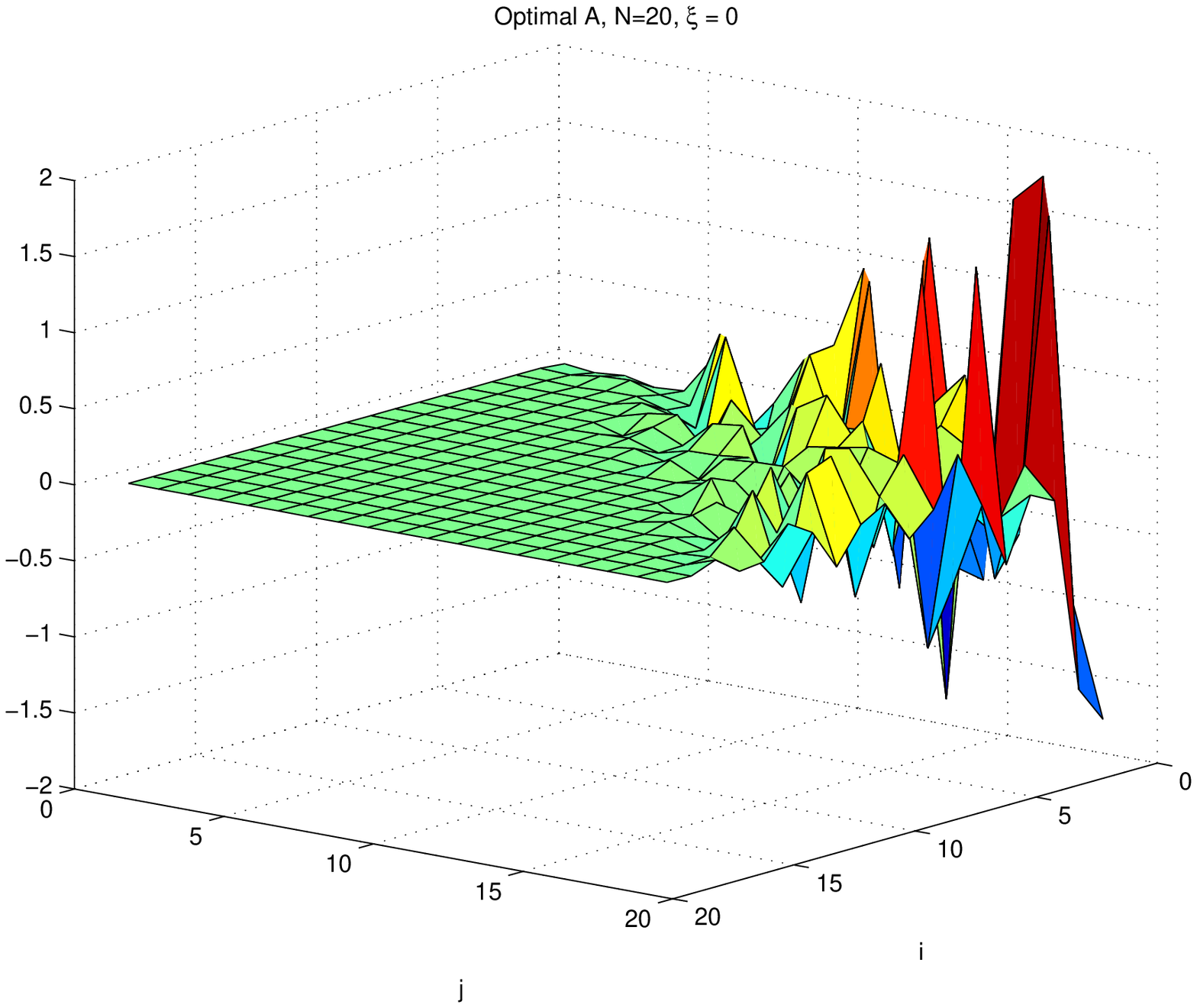}
  \caption{The numerically computed matrix elements $r_{ij}$ and $a_{ij}$ that minimise the
  variance at $\xi = 0$.
  }\label{optimaljagged}
\end{figure}

The optimal elemental coefficients $r_{ij}$ and corresponding
log-spacing coefficients $a_{ij}$ computed for $\xi = 0$ are shown
in Fig.~\ref{optimaljagged}. It can be seen that the coefficients
are small near the $ i \approx j$ diagonal, rising in amplitude near
the top corner $i \approx 1, j \approx N-2$. At all values of $\xi$
investigated, the optimal coefficients had this characteristic.
Moreover all exhibited the decidedly non-smooth character
reminiscent of the measure $\lambda$ in Seger's optimisation
procedure (\cite{segers}).

Although $\xi$-specific optimal combinations have thus been
computed, these are not optimal in any global sense, as their
performance is far from optimal away from the values of $\xi$ at
which they were optimised. This can be seen in Fig.~\ref{optimal},
where the curve D1 shows the performance of the optimal
$\xi=0$ combination falling away rapidly from perfect relative efficiency
for $\xi$ values different from zero.

Fig.~\ref{optimal} also shows the performance of some representative
examples of other linear combinations. Curve A1 is for equal
elemental weights ($r_{ij}=$ constant), showing good performance in
the very heavy-tailed region, but with much lower efficiency for
$\xi$ small or negative.

Curve B1 has $r_{ij} \propto 1/(ii(ii+1))$, and thus gives much
weight to the top row of the matrices, where log-spacings are
measured from the data maximum. This combination is seen to give
excellent performance in the $\xi$ negative region, but has low
efficiency for $\xi$ positive.

Curve C1 (dashed) is for $r_{ij} \propto (jj-ii)^2$. This emulates
the numerically-computed optimals in rising from zero near the
diagonal to larger values in the $i \approx 1, j \approx N-2$ top
corner.  The efficiency is good in the region near $\xi = 1$, but is
low elsewhere, especially for $\xi \approx -1$. Interestingly, the
efficiency stays close to that of the Seger's estimator (+)
throughout.

The relative efficiency of the Seger's estimator is also shown in
Fig.~\ref{optimal}. The reason it can have a relative efficiency
greater than unity (as it does near $\xi \approx 0.5$) is that it is
not strictly within the class of estimators under consideration, in
that it is biased and, moreover, is a nonlinear function of the
log-spacings (since the weights are adaptively selected after an
initial log-spacing-based estimate). It should also be noted that,
despite having the advantages of bias-variance trade-off and access
to a larger class of possibilities, its relative efficiency is
nevertheless comparatively poor for $\xi$ negative.

Finally, it can be seen from Fig.~\ref{optimal} that the linearly
rising combination with $r_{ij} \propto N-j+1$ (shown with circles)
has some plausible claim to being a suitable compromise. It has near
optimal performance in the region $\xi$ from 0 to 3, which is often
of great interest in practice, and although the efficiency falls
somewhat for $\xi$ negative, it does not do so by much. For this
reason, this combination will be considered further throughout the
remainder of the paper.

\subsection{Consistency - preliminary results}

\begin{figure}[h!] \centering
  \includegraphics[width=100mm,keepaspectratio]{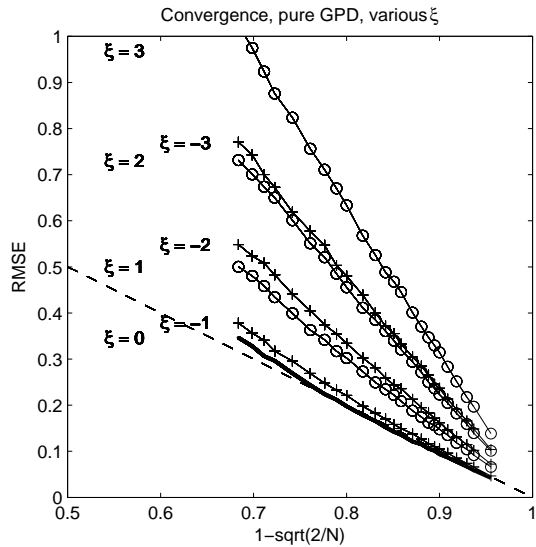}
  \caption{Dependence of root mean square error RMSE on sample size $N$ ($20 \leq N \leq 1000$) for the linearly-rising combination, with samples drawn from GPDs with various tail parameters.
  Positive and negative $\xi$ are indicated with circles and crosses respectively.
  Horizontal axes are scaled to
  bring infinite $N$ to unity, by plotting $1-\sqrt{2/N}$.
  This numerical evidence suggests that, for large $N$, the RMSE decreases in proportion to $1/\sqrt{N}$, with the constant of proportionality
  depending on $\xi$.
  }\label{GPDconsis1}
\end{figure}

Although there is as yet no proof of consistency for any elemental
combination, numerical evidence suggests
that the ``linearly-rising'' combination is consistent for samples
drawn from within the GPD family.
Fig.~\ref{GPDconsis1} shows how the Root Mean Square Error (RMSE)
for the ``linearly-rising'' combination of elementals
decreases as the sample size grows.
At each point on each graph, the RSME was determined numerically from
10,000 samples of size $N$ drawn from a
GPD at various $\xi$ in the range $-3 \leq \xi \leq 3$, with $N$ increasing from 20 to 1000.
At each $\xi$, and for $N$ large, the errors appear to be converging with increasing sample size $N$
at a rate proportional to $1/\sqrt{N}$, with the constant of proportionality dependent on $\xi$.

There is a consistency proof already in existence which has some relevance here, and
covers many elemental combinations for the more general domain of attraction case (which trivially
includes the pure GPD case).
This is Theorem 3.1 of \cite{segers}, which
guarantees weak consistency of many elemental combinations in the $N
\rightarrow \infty$, $k \rightarrow \infty$, $k/N \rightarrow 0$
limit.  To be covered by this theorem, elemental combinations must
be expressible as a mixture of Segers estimators satisfying a condition (Condition 2.5, \cite{segers}), which
re-stated in the notation here, requires the log-spacing
weight matrix $A$ to have zero weights in the vicinity of the
diagonal and of the top row. The linearly-rising combination does not satisfy this condition,
although it is clear that minor adjustments can be made to zero the weights in the appropriate vicinities.

Finally, it could be noted that the presence (or lack) of asymptotic consistency results
is arguably of limited interest if the emphasis, as here, is on providing estimators which perform well
with small or moderately sized samples.

\section{Summary}
`Elemental' absolutely-unbiased, location- and scale-invariant
estimators for the tail parameter $\xi$ of the GPD have been
presented, valid for all $\xi$ and all $N\geq 3$.
The elemental estimators were shown to form a complete basis for
those unbiased, location and scale-invariant estimators of the GPD
tail parameter which are constructed from linear combinations of
log-spacings. Numerical evidence was presented which supports
consistency of at least one elemental combination  for samples
 drawn from the pure GPD family.
 
\bibliography{McRobieGPDbib}

\pagebreak

\input{McRobieGPDappendix1.tex}

\pagebreak

\input{McRobieGPDcompleteness2.tex}

\end{document}

%% file: McRobieGPDappendix1.tex
\section*{Appendix 1: Proof that $\log \tau^{J-1}/t^I$ is unbiased for
the GPD}

Consider a random variable $x$ with distribution function $F(x)$ and
complement $G(x) = 1- F(x)$, and via the probability integral transform,
construct the inverse function $x = u(G)$ which maps an
exceedence probability $G$ to a point $x$ in the data space. Since
$dG = -dF/dx \ dx = -p(x) \ dx$ where $p(x)$ is the density, the
expected value of any function $h(x)$ may be evaluated by transforming integrals
over $x$ to integrals over $G$, using :
\begin{equation}
\langle h(x) \rangle = \int_{\forall x} h(x) \ p(x) \ dx = \int_0^1
h(u(G)) \ dG \label{A1}
\end{equation}
A consequence of the well-known uniform density of the $G$'s is that
the second integral (over $G$) can be considerably simpler than the
first integral (over $x$). For the GPD, integrals via $p(x)$ lead to
lengthy expressions involving hypergeometric (Lauricella) functions.
Although a proof that the elemental estimators are unbiased has been
constructed by that route, the approach via the transformed
integrals over $G$ is considerably simpler, and is presented here.

Consider an ordered sample $\bX$ of $n$ data points drawn from the
distribution $F$, ordered such that $X_n \leq X_{n-1} \leq \ldots
\leq X_2
  \leq X_1$.  The expected value of any function $h(\bX)$ is
\begin{equation}
\langle h(\bX) \rangle =  n! \ \int_0^1 dG_n \ldots \int_0^{G_2}
dG_1 \ h(u(\bG))
\end{equation}
the integral being over the $n$-dimensional unit simplex containing
all possible $G$.

\begin{eqnarray}
\text{The GPD has distribution function} \ F(x) & = & 1 - \left( 1+
\frac{\xi(x-\mu)}{\sigma} \right)^{-1/\xi} \\
 \text{such that} \ x & = &
u(G) = \mu + \frac{\sigma}{\xi}(G^{-\xi} - 1)
\end{eqnarray}
Depending whether $\xi$ is positive or negative, the expected value
of the log-spacing between the $i$th and $j$th order statistics is
\begin{equation} \langle \log( X_i - X_j)
\rangle =
\begin{cases}
\langle \log (G_j^\gamma - G_i^\gamma) \rangle - \gamma \langle \log
G_i G_j \rangle + \log \frac{\sigma}{\gamma} & \text{for} \ \xi = \gamma, \ \ \gamma >0 \\
 \langle \log (G_j^\gamma -
G_i^\gamma) \rangle \hspace{2.1cm} + \log \frac{\sigma}{\gamma} &
\text{for} \ \xi = -\gamma, \ \gamma >0
\end{cases}
\label{squigneg}
\end{equation}

Consider an estimator $ \ \hat{\xi}(\bX) = \sum_{i,j} a_{ij}
\log(X_{i} - X_{j})$, a linear combination of log-spacings.
For $\hat{\xi}$ to be scale invariant, the weights $a_{ij}$ must sum
to zero to remove the $\sigma$ dependence in Eqn.~(\ref{squigneg}).
Moreover, for $\hat{\xi}$ to be unbiased for both positive and
negative $\xi$, Eqn.~(\ref{squigneg}) requires
\begin{equation}
\sum_{i,j} a_{ij} \langle \log ( G_j^\gamma - G_i^\gamma ) \rangle =
-\gamma \ \ \ \text{and} \ \ \ - \sum_{i,j} a_{ij} \gamma \langle
\log G_i G_j \rangle = 2\gamma
\end{equation}
To determine the expected value of any function $h(X_i, X_j)$, all
other $G$ variables can be integrated out, leaving 
\begin{equation}
\langle h(X_i,X_j) \rangle = C_{ij} \int_0^1 \ dG_j \int_0^{G_j}
dG_i \  \ \ G_i^{i-1}(1-G_j)^{n-j}(G_j-G_i)^{j-i-1} \ h(u(G_i),
u(G_j))\label{equation8}
\end{equation}
where $C_{ij} = n!/((i-1)! (n-j)!(j-i-1)!)$. For example
\begin{equation}
\sum a_{ij} \langle \log G_j \rangle  = \sum a_{ij} C_{ij} .
\int_0^1 G_j^{j-1}(1-G_j)^{n-j} \log G_j \ dG_j \ . \ \int_0^1
\phi^{i-1}(1-\phi)^{j-i-1} \ d\phi \label{oldeqn9}
\end{equation}
where $\phi = G_i/G_j$. The $\phi$ integral leads immediately to the
beta function $B(i, j-i)$. For the $G_j$ integral, standard Mellin
transform theory gives
\begin{eqnarray}
\int_0^1 G_j^{j-1}(1-G_j)^{n-j} \log G_j \ dG_j  & = & \left[
\frac{d}{ds}\int_0^1 G_j^{s-1+j-1}(1-G_j)^{n-j} \ dG_j \right]_{s=1}
\nonumber \\ & = &   B(j, n-j+1) \left( \psi(j) - \psi(n+1) \right)
\label{oldeqn10}
\end{eqnarray}
where $\psi$ is the digamma function, the derivative of the
logarithm of the gamma function. Both beta functions have integer
arguments and may be expressed as ratios of factorials. These cancel
with the leading factorial terms $C_{ij}$, such that the expected
value of the weighted sum is
\begin{equation}
\sum a_{ij} \langle \log G_j \rangle
 = \sum a_{ij} \left( \psi(j) - \psi(n+1) \right)
 =  \sum a_{ij} \psi(j) \label{weightedsum}
\end{equation}
the last step following from $\sum a_{ij} = 0$. It follows similarly
that
\begin{equation}
\sum a_{ij} \langle \log G_i \rangle
 = \sum a_{ij} \left( \psi(i) - \psi(n+1) \right)
 =  \sum a_{ij} \psi(i) \label{weightedsum1}
\end{equation}

For an individual elemental $\hat{\xi}_{IJ}$ the weights $a_{ij}$
and the relevant values of $i$ and $j$ are given in Table
\ref{tableij}.

\begin{table}[h!]
\caption{The indices of an elemental estimator, and, in the final
column, the weights. \label{tableij}} \centering
\begin{tabular}{|c||c|c||c||}
  \hline
  term                      & $i$ &     $j$ &  $a_{ij}$ \\ \hline
   $\log{(X_I - X_{J-1})} $ & $I$ &   $J-1$ &  $J-1$ \\
   $\log{(X_I - X_{J})}$    & $I$ &     $J$ &  $-(J-I-1)$ \\
   $\log{(X_{I+1} - X_J )}$ & $I+1$ &   $J$ &  $-I$ \\
  \hline
\end{tabular}
\end{table}

Using the
relation $\psi(1+x) = \psi(x)+ 1/x$ we obtain the contributions from
the individual elemental $\hat{\xi}_{IJ}$ to be
\begin{eqnarray}
\sum a_{ij} \psi(j)  & = & (J-1)\psi(J-1) -(J-I-1) \psi(J) -I
\psi(J)  \nonumber \\ & = &  (J-1)( \psi(J-1) -  \psi(J))  =  -1 \\
\sum a_{ij} \psi(i)  & = & (J-1)\psi(I) -(J-I-1) \psi(I) -I
\psi(I+1) \nonumber \\ & = & \ I( \psi(I) -  \psi(I+1)) \ = \ -1
\label{apsij}
\end{eqnarray}

These are exactly what is required to show that any unit-sum linear
combination $\hat{\xi}$ of elementals provides an unbiased estimate
for $\xi$ via the $\langle \log G_i^{\gamma} \rangle$ and $\langle
\log G_j^{\gamma} \rangle$ terms. Explicitly, with $\gamma = | \xi
|$, we have
\begin{eqnarray}
\langle \hat{\xi} \rangle  & = & \left\{  \begin{array}{ll}  \
\gamma \sum a_{ij} \langle \log G_j \rangle + \sum a_{ij} \langle
\log( 1-
\phi^\gamma) \rangle & \mbox{for $\xi < 0 $ } \\
- \gamma \sum a_{ij} \langle \log G_i \rangle + \sum a_{ij} \langle
\log( 1- \phi^\gamma) \rangle & \mbox{for $\xi > 0$ } \end{array} \right. \\
& = & \xi + \sum a_{ij} \langle \log( 1- \phi^\gamma) \rangle  \ \ \
\end{eqnarray}
It only remains to prove that the second term is zero. This term
involves somewhat more complicated integrals leading to sums of
digamma functions with non-integer arguments dependent on $\gamma$.
Explicit evaluation of these can, however, be avoided here by
observing how the elemental terms combine and cancel.

From Eqn.~(\ref{equation8}), we obtain
\begin{eqnarray}
\langle \log( 1- \phi_{ij}^\gamma) \rangle    &  = & C_{ij} D_{j} B_{ij} \nonumber \\
\text{where} \ \ D_{j} & = & \int_0^1 y^{j-1} (1-y)^{n-j} \ dy  =  \frac{(j-1)!(n-j)!}{n!}  \nonumber \\
\text{and} \ \ \ B_{ij} & = &  \int_0^1
\phi^{i-1}(1-\phi)^{j-i-1}  \log(1-\phi^\gamma) \ d \phi
\end{eqnarray}

Summing over a single elemental using the weights $a_{ij}$ given in Table \ref{tableij}, and collecting the leading $a_{ij} C_{ij}Dj$ terms gives
\begin{equation}
\sum a_{ij} \langle \log( 1- \phi^\gamma) \rangle  
  =   \frac{(J-1)!}{(I-1)!(J-I-2)!}  (B_{I,J-1} - B_{I,J}  - B_{I+1,J}) \\
\end{equation}
Now
\begin{eqnarray}
\lefteqn{
  B_{I,J-1} - B_{I,J}  - B_{I+1,J}    =  }  \\
 & &  \int_0^1 \lbrace
\phi^{I-1}(1-\phi)^{J-I-2} - \phi^{I-1}(1-\phi)^{J-I-1} - \phi^I
(1-\phi)^{J-I-2} \rbrace \log(1-\phi^\gamma) \ d \phi \nonumber
\end{eqnarray}
where the integrand contains the factor $ \left[ 1 - (1-\phi) - \phi
\right] = 0$,  thus
\begin{equation} \sum a_{ij} \langle \log( 1-
\phi^\gamma) \rangle = 0 \ \ \ \ \ \ \ \forall \xi \neq 0
\end{equation}
This completes the proof for all nonzero $\xi$.

The proof for $\xi = 0 $ follows similarly, using $G(x) = \exp{(-(x-\mu)/\sigma)}$, whence, using
 Eqn.~(\ref{equation8}), we obtain
 \begin{equation}
 \sum a_{ij} \langle \log (X_i-X_j) \rangle = \sum a_{ij} \langle \log \sigma \rangle + \sum a_{ij} \langle \log [-\log \phi_{ij} ] \rangle
\end{equation}
The first term is trivially zero due to the zero sum of the elemental $a_{ij}$, and the second term is zero for the same reason that
$\sum a_{ij} \langle \log (1 - \phi^\gamma) \rangle$ is zero, as shown above, in that the elementals combine terms in such a way as to eliminate
the expectation of any function of the $\phi_{ij}$.

This
completes the proof that for any elemental $\hat{\xi_{IJ}}$ , the expectation $\langle \hat{\xi_{IJ}} \rangle = \xi$ for all $\xi$.

%% file: McRobieGPDcompleteness2.tex
\section*{Appendix 2: Completeness}

Here we prove that for an estimator $\hat{\xi}$ of the tail
parameter $\xi$ of the GPD, with the preconditions that $\hat{\xi}$
is i) a linear combination of log-spacings, ii) absolutely-unbiased
for all $\xi$ and iii) location- and scale-invariant,  then
$\hat{\xi}$ may be expressed as a linear combination of the
elementals. That is, the elementals form a complete basis for the
set of invariant unbiased log-spacing estimators of the GPD tail
parameter.

Precondition ii requires that $\hat{\xi}$ must be unbiased at both
$\xi = \gamma$ and $\xi = -\gamma$ for any $\gamma>0$. This symmetry
is embodied in Eqn.~(\ref{squigneg}), addition and subtraction of
which (together with precondition iii and elementary integrations
such as Eqns~\ref{oldeqn9}-\ref{weightedsum}) leads to the
requirements on the log-spacing weights $a_{ij}$ that
\begin{eqnarray}
\lefteqn { \gamma \sum a_{ij} ( \langle \log G_j \rangle + \langle
\log G_i \rangle ) } \nonumber \\ &
= & \gamma \sum a_{ij} \left[ \psi (j) + \psi(i)\right]  = -2 \gamma \\
\lefteqn{ \gamma \sum a_{ij} ( \langle \log G_j \rangle - \langle
\log G_i \rangle ) } \nonumber \\ & = & \gamma \sum a_{ij} \left[
\psi (j) - \psi(i)\right] = - 2 \sum a_{ij} \langle \log
(1-\phi_{ij}^\gamma) \rangle \label{mydiff}
\end{eqnarray}
where $\psi$ is the digamma function and $\sum$ means sum over $i$
from 1 to $N-1$ and $j$ from $i+1$ to $N$.
\\
We now prove that the preconditions imply that the right-hand side
of Eqn~(\ref{mydiff}) is zero. Eqn~(\ref{mydiff}) may be written as
\begin{equation}
c_1 \gamma = \sum a_{ij} \langle \log(1-\phi_{ij}^\gamma) \rangle :=
I_1  \ \text{  where  } c_1 = -(1/2)\sum a_{ij} \left[ \psi(j) -
\psi(i) \right] \label{Idefintion}
\end{equation}
Similar to Eqn~(\ref{oldeqn9}), each term in the $I_1$ summation may
be considered individually and all variables irrelevant to that term
may be integrated out to give
\begin{equation}
I_1 = \sum_{ij} a_{ij} C_{ij} \int_0^1 G_j^{j-1} (1-G_j)^{n-j} \
dG_j \ . \ \int_0^1 {\phi_{ij}}^{i-1}(1-\phi_{ij})^{j-i-1}
\log(1-\phi_{ij}^\gamma) \ d\phi_{ij} \label{myI1}
\end{equation}
The $G_j$ integral gives a beta function $B(j, n-j+1)$ which
combines with the $C_{ij}$ term to give a factor $1/B(i,j-i)$. Since
each expectation integral over $\phi_{ij}$ is definite, we can 
set each $\phi_{ij} = \phi$. 
Each product $\phi^{i-1}(1-\phi)^{j-i-1}$ involves integer exponents
and can thus be expressed as a polynomial of degree $j-2$. Passing
the summation through the integral sign, the various polynomials can
be collected into a single polynomial $p_{n-2}(\phi) =
\sum_{k=0}^{n-2} b_k \phi^k $ of degree $n-2$. Passing the summation
back out of the integral gives:
\begin{equation}
I_1 = \sum_{k=0}^{n-2} b_k g_k(\gamma) \ \text{  where } \
g_k(\gamma): = \int_0^1 \phi^k \log(1-\phi^\gamma) \ d\phi
\label{gkdef}
\end{equation}
\\
The binomial expansion of each $(1-\phi)^{j-i-1}$ factor in
Eqn.~\ref{myI1} leads to the total polynomial
\begin{equation}
p_{n-2}(\phi) = \sum_{i=1}^{n-1} \sum_{j=i+1}^{n}
\frac{(j-1)!}{(i-1)!(j-i-1)!} a_{ij} \sum_{q=0}^{j-i-1} \frac{(-1)^q
(j-i-1)!}{q!(j-i-1-q)!} \phi^{(i-1)+q}
\end{equation}
Collecting together equal powers of $\phi$ 
gives
\begin{equation}
p_{n-2}(\phi) = \sum_{k=0}^{n-2} \phi^k \sum_{i=1}^{k+1}
\sum_{j=k+2}^n \frac{(j-1)!}{(i-1)!} a_{ij}
\frac{(-1)^{k-i-1}}{(k-i+1)!(j-k-2)!}
\end{equation}
such that the polynomial coefficients $b_k$ are
\begin{eqnarray}
b_k  & = & \sum_{i=1}^{k+1} \sum_{j=k+2}^n \frac{(j-1)!}{(i-1)!}
a_{ij} \frac{(-1)^{k-i-1}}{(k-i+1)!(j-k-2)!} \nonumber \\
&  = & (k+1) \sum_{j=k+2}^n \left(
\begin{array}{c} j-1 \\ k+1 \end{array} \right)
 \sum_{i=1}^{k+1} \left(
 \begin{array}{c} k \\ i-1 \end{array} \right) (-1)^{k-i-1} a_{ij}
 \label{bkweights}
\end{eqnarray}
Having determined the polynomial coefficients, we consider the
integrals $g_k(\gamma)$ of the various $\phi^k$ terms, as defined in
Equation (\ref{gkdef}). The transformation $\phi^\gamma = (1-\rho)$
leads to
\begin{eqnarray}
g_k(\gamma) & = & \frac{1}{\gamma} \int_0^1
(1-\rho)^{\frac{k+1}{\gamma} -1} \log \rho \ d\rho =
\frac{1}{\gamma} \left[ \frac{d \ }{ds} \int_0^1
\rho^{s-1}(1-\rho)^{\frac{k+1}{\gamma}-1}  \ d\rho \right]_{s=1} \nonumber  \\
& = & \frac{1}{\gamma} B(1,\frac{k+1}{\gamma})
\left[\psi(1)-\psi(1+\frac{k+1}{\gamma})\right] \nonumber \\ & = &
\frac{1}{k+1} \left[\psi(1)-\psi(1+\frac{k+1}{\gamma})\right]
\end{eqnarray}

We thus have
\begin{equation}
I_1 = c_1 \gamma = \sum_{k=0}^{n-2} b_k . g_k(\gamma) =
\sum_{k=0}^{n-2}
\frac{b_k}{k+1}\left[\psi(1)-\psi(1+\frac{k+1}{\gamma})\right]
\end{equation}
\\
If this is true at all $\gamma > 0$, then it must be true for
$\gamma$ large, ($\gamma = 1/\epsilon$, $\epsilon > 0$, $\epsilon$
small). There
\begin{equation}
g_k \left( \gamma = \frac{1}{\epsilon} \right) = \frac{1}{k+1}
\left[ \psi (1)-\psi (1+(k+1) \epsilon) \right] \nonumber
\end{equation}
Since the digamma function is well-behaved (i.e. infinitely
differentiable) near $\psi(1)$, we may take the Taylor series
$\psi(1+\delta) = \psi(1) + \psi '(1) \delta + O(\delta^2)$ to
obtain
\begin{equation}
g_k(\gamma) =  -\psi ' (1)  \epsilon + O(\epsilon^2)
\end{equation}
thus
\begin{equation}
I_1 = c_1 \gamma = \frac{c_1}{\epsilon} = \sum_{k=0}^{n-2} b_k
g_k(\frac{1}{\epsilon}) = -\epsilon \psi ' (1) \sum_{k=0}^{n-2}
{b_k} + O(\epsilon^2)
\end{equation}
whence
\begin{equation}
c_1 = - \epsilon^2 \psi ' (1) \sum_{k=0}^{n-2} {b_k} + O(\epsilon^3)
\label{c1condition}
\end{equation}
Since $c_1$ is a constant wrt $\gamma$ and since (nonzero) $\epsilon
= 1/\gamma$ can be arbitrarily small,
 we thus infer from \ref{c1condition} that the preconditions imply that $c_1 = 0$.

Finally, since the preconditions imply that $I_1 =  c_1 \gamma = 0$
and $I_1 =\sum_{k=0}^{n-2} b_k \ g_k(\gamma)$
the independence of the $N-1$ functions $g_k(\gamma)$ implies that
\begin{equation}
b_k = 0, \ \ \forall k \in \lbrace 0,1,\ldots, n-2 \rbrace
\end{equation}
These are the $N-1$ constraints needed to reduce the dimension of
the problem down to that spanned by the elementals.

That the elementals are contained within this subspace can be
readily checked by gathering, within each $b_k$ summation, the terms
associated with each elemental crosshair of the grid $G$. Each $b_k
= 0$ constraint corresponds to a weighted summation over a
subrectangle of the $A$ matrix as illustrated in
Table~\ref{constraintmatrix}.

\begin{table}[h!]
\centering
\begin{tabular}{ |ccccccc|}
  \hline
  . & $a_{12}$ & $a_{13}$ & $a_{14}$ & \cs $a_{15}$ & \cs $a_{16}$ & \cs $a_{17}$ \\
    &          & $a_{23}$ & $a_{24}$ & \cs $a_{25}$ & \cs $a_{26}$ & \cs $a_{27}$ \\
    &          &          & $a_{34}$ & \cs $a_{35}$ & \cs $a_{36}$ & \cs $a_{37}$ \\
    &          &          &          & \cs $a_{45}$ & \cs $a_{46}$ & \cs $a_{47}$ \\
    &          &          &          &              &     $a_{56}$ &     $a_{57}$ \\
    &          &          &          &              &              &     $a_{67}$ \\
    &          &          &          &              &              &              \\
  \hline
\end{tabular}
\caption{Each constraint $b_k = 0$ is a summation over a rectangle
of the $A$ matrix. For $N=7$, the rectangle of terms for $k=3$ is
shown.} \label{constraintmatrix}
\end{table}

Consider an element $a_{IJ}$ away from the rectangle boundaries
(such as $a_{26}$ in Table~\ref{constraintmatrix}).  Consider the
$b_k$ summation terms associated with the elemental
$\hat{\xi}_{IJ}$, which thus involves the term
 $a_{IJ}$, the term $a_{I,J-1}$ to its left and the term $a_{I+1,J}$
below it. The $b_k$ summation weights given to each of these terms
can be obtained from Eqn~\ref{bkweights} as
\begin{eqnarray}
b_{k,(I,J)}                & = & \frac{(J-1)! (+1)}{(I-1)!(k-I+1)!(J-k-2)!}\\
b_{k,(I,J-1)}              & = & \frac{(J-2)! (+1)}{(I-1)!(k-I+1)!(J-k-3)!}\\
\text{and} \ b_{k,(I+1,J)} & = &
\frac{(J-1)!(-1)}{(I)!(k-I)!(J-k-2)!}
\end{eqnarray}

 The elemental $\hat{\xi}_{IJ}$
contributes in proportions $-(J-I-1)$, $(J-1)$ and $-I$ to $a_{IJ}$,
$a_{I,J-1}$ and $a_{I+1,J}$ respectively. The weights in the $b_k$
summation are thus such as to eliminate the contribution from the
elemental $\xi_{IJ}$, since it readily follows from the above that
\begin{eqnarray}
\lefteqn{-(J-I-1)b_{k,(I,J)}+ (J-1)b_{k,(I,J-1)}+ (-I)b_{k,(I+1,J)}
} \\
& \propto &
-(J-I-1)+(k-I+1)+(J-k-2) \nonumber\\
& = & 0
\end{eqnarray}

(The proof for elements near the boundaries of the $b_k$ rectangle
is similar).

Since any elemental satisfies the constraints then so does any
linear combination thereof, and since there are $(N-1)(N-2)/2$
elementals and they are independent, it follows that they form a
complete basis for those estimators of the GPD tail index that
satisfy the preconditions given.

%% file: McRobieGPD14April2013arxiv1.bbl
\begin{thebibliography}{13}
\expandafter\ifx\csname natexlab\endcsname\relax\def\natexlab#1{#1}\fi
\expandafter\ifx\csname url\endcsname\relax
  \def\url#1{\texttt{#1}}\fi
\expandafter\ifx\csname urlprefix\endcsname\relax\def\urlprefix{URL }\fi

\bibitem[{Castillo and Daoudi(2009)}]{castillo}
Castillo, J., Daoudi, J., 2009. Estimation of the generalized {P}areto
  distribution. Statistics and Probability Letters 79, 684--688.

\bibitem[{Coles(2001)}]{coles}
Coles, S., 2001. An Introduction to Statistical Modelling of Extreme Values.
  Springer, London.

\bibitem[{de~Zea~Bermudez and Kotz(2010)}]{bermudez1}
de~Zea~Bermudez, P., Kotz, S., 2010. Parameter estimation for the generalized
  {P}areto distribution - parts i and ii. J. Stat. Planning and Inference
  140~(6), 1353--1388.

\bibitem[{Dekkers et~al.(1989)Dekkers, Einmahl, and deHaan}]{Dekkers}
Dekkers, A. L.~M., Einmahl, J. H.~J., deHaan, L., 1989. A moment estimator for
  the index of an extreme value distribution. Annals of Statistics 17,
  1833--1855.

\bibitem[{Drees(1998)}]{drees2}
Drees, H., 1998. A general class of estimators of the extreme value index. J.
  Stat. Planning and Inference 66, 95--112.

\bibitem[{Embrechts et~al.(1999)Embrechts, Kl\"{u}ppelberg, and
  Mikosch}]{embrechts}
Embrechts, P., Kl\"{u}ppelberg, C., Mikosch, T., 1999. Modelling Extreme Events
  for Insurance and Finance. Springer, Berlin.

\bibitem[{Hill(1975)}]{hill}
Hill, B.~M., 1975. A simple general approach to inference about the tail of a
  distribution. Annals of Statistics 3, 1163--1174.

\bibitem[{Hosking and Wallis(1987)}]{hosking}
Hosking, J. R.~M., Wallis, J.~R., 1987. Parameter and quantile estimation for
  the {G}eneralized {P}areto {D}istribution. Technometrics 29~(3), 339--349.

\bibitem[{Pickands(1975)}]{pickands}
Pickands, J., 1975. Statistical inference using extreme order statistics. Ann.
  Statist. 3, 119--131.

\bibitem[{Reiss and Thomas(2001)}]{reiss}
Reiss, R.-D., Thomas, M., 2001. Statistical Analysis of Extreme Values: with
  applications to insurance, finance, hydrology and other fields.
  Birkh\"{a}user, Basel.

\bibitem[{Segers(2005)}]{segers}
Segers, J., 2005. Generalized {P}ickands estimators for the extreme value
  index. J. Stat. Planning and Inference 128~(2), 381--396.

\bibitem[{Smith(1987)}]{smith}
Smith, R.~L., 1987. Estimating tails of probability distributions. Annals of
  Statistics 15~(3), 1174--1207.

\bibitem[{Yun(2002)}]{yun}
Yun, S., 2002. On a generalized {P}ickands estimator of the extreme value
  index. J. Stat. Planning and Inference 102, 389--409.

\end{thebibliography}
